\documentclass[12pt,thmsa]{article}
\usepackage{amssymb}
\begin{document}

\author{David Carf\`{i}, Daniele Schilir\`{o}}
\title{A model of coopetitive games and the Greek crisis}
\date{}
\maketitle

\begin{abstract}
In the present work we propose an original analytical model of coopetitive
game. We shall try to apply this analytical model of coopetition - based on normal form game
theory and conceived at a macro level - to the Greek crisis, suggesting
feasible solutions in a super-cooperative perspective for the divergent interests
which drive the economic policies in the euro area.
\end{abstract}

\bigskip


\bigskip

\bigskip

\section{\textbf{Introduction}}

\bigskip

In this contribution we focus on the Greek crisis, since we know that Greece
is still in a very difficult economic situation, due to its lack of
competitiveness and is at risk of insolvency, because of its public finance
mismanagement. Although the EU Governments and IMF have recently approved
more substantial financial aids to cover the refinancing needs of Greece
until 2014, in exchange of a serious and tough austerity program. Germany,
on the other hand, is the most competitive economy of the Euro Area and has
a large trade surplus with Greece and other Euro partners; hence significant
trade imbalances occur within the Euro Area.

\bigskip

The main purpose of our paper is to explore win-win solutions for Greece and
Germany, involving a German increasing demand of a pre-determined Greek exports. We do not
analyze the causes of the financial crisis in Greece and its relevant
political and institutional effects on the European Monetary Union. Rather
we concentrate on stability and growth, which should drive the economic
policy of Greece and the other Euro countries.

\bigskip

\textbf{Organization of the paper.} The work is organized as follows:

\begin{itemize}
\item  section 2 examines the Greek crisis, suggesting a possible way out to
reduce the intra-eurozone imbalances through coopetitive solutions within a
growth path;

\item  sections from 3 to 6 provides an original model of coopetitive game
applied to the Eurozone context, showing the possible coopetitive solutions;

\item  conclusions end up the paper.
\end{itemize}

Introduction and Section 2 of this paper are written by D. Schilir\`{o},
sections from 3 to 6 are written by D. Carf\`{i}; conclusions are written
by the two authors.

\bigskip

\textbf{Acknowledgments.} We wish to thank Daniela Baglieri, Albert E.
Steenge and three anonymous referees for their helpful comments and
suggestions.

\bigskip

\section{\textbf{The Greek Crisis and the coopetitive solution}}

\bigskip

The deep financial crisis of Greece, which was almost causing the default of
its sovereign debt, has revealed the weaknesses of Greek economy,
particularly its lack of competitiveness, but also the mismanagement of the
public finance and the difficulties of the banking sector.

\bigskip

\subsection{The crisis and the Greek economy}

\bigskip

With the outbreak of the global crisis of 2008-2009, Greece relied on state
spending to drive growth. Moreover, the country has accumulated a huge
public debt (over 320 billion euros in 2010). This has created deep concerns
about its fiscal sustainability and its financial exposition has prevented
the Greek government to find capitals in the financial markets. In addition,
Greece has lost competitiveness since joining the European Monetary Union
and, because of that, Greek's unit labor cost rose 34 percent from 2000 to
2009. The austerity measures implemented by the Greek government, although
insufficient, are hitting hard the Greek economy, since its growth is
expected to be negative also this year (2011), making the financial recovery
very problematic [Mussa, 2010]. Furthermore, Greece exports are much less
than imports, so the trade balance shows a deficit around 10\%. Therefore,
the focus of economic policy of Greece should become its productive system
and growth must be the major goal for the Greek economy in a medium term
perspective. This surely would help its re-equilibrium process.

\bigskip

\subsection{The soundest European economy: Germany}

\bigskip

Germany, on the other hand, is considered the soundest European economy. It
is the world's second-biggest exporter, but its wide commercial surplus is
originated mainly by the exports in the Euro area, that accounts for about
two thirds. Furthermore, since 2000 its export share has gradually increased
vis-\`{a}-vis industrial countries. Thus Germany's growth path has been
driven by exports. We do not discuss in this work the factors explaining
Germany's increase in export share, but we observe that its international
competitiveness has been improving, with the unit labor cost which has been
kept fairly constant, since wages have essentially kept pace with
productivity. Therefore the prices of the German products have been
relatively cheap, favoring the export of German goods towards the euro
countries and towards the markets around the world, especially those of the
emerging economies (China, India, Brasil, Russia). Finally, since 2010
Germany has recovered very well from the 2008-2009 global crisis and it is
growing at a higher rate than the others Euro partners.

\bigskip

Therefore we share the view that Germany (and the other surplus countries of
the Euro area) should contribute to overcome the crisis of Greek economy
stimulating its demand of goods from Greece and relying less on exports
towards the Euro area in general. Germany, as some economists as Posen
[2010] and Abadi [2010] underlined, has benefited from being the anchor
economy for the Eurozone over the last 11 years. For instance, in 2009,
during a time of global contraction, Germany has been a beneficiary, being
able to run a sustained trade surplus with its European neighbors. Germany
exported, in particular, 6.7 billions euro worth of goods to Greece, but
imported only 1.8 billion euro worth in return.

\bigskip

\subsection{A win-win solution for Greece and Germany}

\bigskip

Thus we believe that an economic policy that aims at adjusting government
budget and trade imbalances and looks at improving the growth path of the
real economy in the medium and long term in Greece is the only possible one
to assure a stable re-balancing of the Greek economy and also to contribute
to the stability of the whole euro area [Schilir\`{o}, 2011]. As we have
already argued, German modest wage increases and weak domestic demand
favored the export of German goods towards the euro countries. We suggest,
in accordance with Posen [2010], a win-win solution (a win-win solution is
the outcome of a game which is designed in a way that all participants can
profit from it in one way or the other), which entails that Germany, which
still represents the leading economy, should re-balance its trade surplus
and thus ease the pressure on the southern countries of the euro area,
particularly Greece. Obviously, we are aware that this is a mere hypothesis
and that our framework of coopetition is a normative model. However, we
believe that a cooperative attitude must be taken within the members of the
European monetary union. Thus we pursue our hypothesis and suggest a model
of coopetitive game as an innovative instrument to analyze possible
solutions to obtain a win-win outcome for Greece and Germany, which would
also help the whole economy of the euro area.

\bigskip

\subsection{Our coopetitive model}

\bigskip

The two strategic variables of our model are investments and exports for
Greece, since this country must concentrate on them to improve the structure
of production and its competitiveness, but also shift its aggregate demand
towards a higher growth path in the medium term. Thus Greece should focus on
innovative investments, specially investments in knowledge [Schilir\`{o},
2010], to change and improve its production structure and to increase its
production capacity and its productivity. As a result of that its
competitiveness will improve. An economic policy that focuses on investments
and exports, instead of consumptions, will address Greece towards a
sustainable growth and, consequently, its financial reputation and economic
stability will also increase. On the other hand, the strategic variable of
our model for Germany private consumption and imports.

\bigskip

The idea which is driving our model to solve the Greek crisis is based on a
notion of coopetition where the cooperative aspect will prevail. Thus we are
not talking about a situation in which Germany and Greece are competing in
the same European market for the same products, rather we are assuming a
situation in which Germany stimulates its domestic demand and, in doing so,
will create a larger market for products from abroad. We are also envisaging
the case where Germany purchases a greater quantity of Greek products, in
this case Greece increases its exports, selling more products to Germany.
The final results will be that Greece will find itself in a better position,
but also Germany will get an economic advantage determined by the higher
growth in the two countries. In addition, there is the important advantage
of a greater stability within the European Monetary system. Finally our
model will provide a new set of tools based on the notion of coopetition,
that could be fruitful for the setting of the euro area economic policy
issues.

\bigskip

\subsection{The coopetition in our model}

\bigskip

The concept of coopetition was essentially devised at micro-economic level
for strategic management solutions by Brandenburger and Nalebuff [1995], who
suggest, given the competitive paradigm [Porter, 1985], to consider also a
cooperative behavior to achieve a win-win outcome for both players.
Therefore, in our model, coopetition represents the synthesis between the
competitive paradigm [Porter, 1985] and the cooperative paradigm [Gulati,
Nohria, Zaheer, 2000; Stiles, 2001]. Coopetition is, in our approach, a
complex theoretical construct and it is the result of the interplay between
competition and cooperation. Thus, we suggest a model of coopetitive games,
applied at a macroeconomic level, which intends to offer possible solutions
to the partially divergent interests of Germany and Greece in a perspective
of a cooperative attitude that should drive their policies.

\bigskip

\section{\textbf{Coopetitive games}}

\bigskip

\subsection{\textbf{Introduction}}

\bigskip

In this paper we develop and apply the mathematical model of \emph{%
coopetitive game} introduced by David Carf\`{i} in \cite{ca-sch1} and  \cite{ca-Sing6}. The idea of
coopetitive game is already used, in a mostly intuitive and non-formalized
way, in Strategic Management Studies (see for example Brandenburgher and Nalebuff).

\bigskip

\textbf{The idea.} A coopetitive game is a game in which two or more players (participants)
can interact \emph{cooperatively and non-cooperatively at the same time}.
But even Brandenburger and Nalebuff, creators of coopetition, did not
defined precisely a \emph{quantitative way to implement coopetition} in Game
Theory context.

\bigskip

The problem in Game Theory to implement the notion of coopetition is:

\begin{itemize}
\item  \emph{how do, in normal form games, cooperative and non-cooperative
interactions can live together simultaneously, in a Brandenburger-Nalebuff
sense?}
\end{itemize}

\bigskip

Indeed, consider a classic two-player normal-form gain game $G=(f,>)$ - such
a game is a pair in which $f$ is a vector valued function defined on a
Cartesian product $E\times F$ with values in the Euclidean plane $\Bbb{R}^{2} $ and $>$ is the natural strict sup-order of the Euclidean plane itself. Let $E$ and $F$ be the strategy sets of the two players in the game $G$. The two players
can choose the respective strategies $x\in E$ and $y\in F$ cooperatively
(exchanging information) or not-cooperatively (not exchanging informations),
but these two ways are mutually exclusive, at least in normal-form games: the two ways
cannot be adopted simultaneously  in the model of normal-form
game (without using convex probability mixtures, but this is not the
way suggested by Brandenburger and Nalebuff in their approach). There is no room, in the classic normal game model, for a simultaneous
(non-probabilistic) employment of the two behavioral extremes \emph{cooperation} and 
\emph{non-cooperation}.

\bigskip

\textbf{Towards a possible solution.} David Carf\`{i} (\cite{ca-sch1} and  \cite{ca-Sing6}) has
proposed a manner to pass this \emph{impasse}, according to the idea of
coopetition in the sense of Brandenburger and Nalebuff:

\begin{itemize}
\item  in a Carf\`{i}'s coopetitive game model, the players of the game have their
respective strategy-sets (in which they can choose cooperatively or not
cooperatively) and a common strategy set $C$ containing other strategies
(possibly of different type with respect to those in the respective classic strategy sets) that \emph{must be chosen cooperatively}. This strategy set $C$ can also be structured
as a Cartesian product (similarly to the profile strategy space of normal form games), but in any case the strategies belonging to this new
set $C$ must be chosen cooperatively.
\end{itemize}

\bigskip

\subsection{\textbf{The model for }$n$\textbf{\ players}}

\bigskip

We give in the following the definition of coopetitive game proposed by
Carf\`{i} (in \cite{ca-sch1} and  \cite{ca-Sing6}).

\bigskip

\textbf{Definition (of }$n$\textbf{-player coopetitive game).}\emph{\ Let }$E=(E_{i})_{i=1}^{n} $ \emph{be a finite }$n$\emph{-family of non-empty sets and let }$C$\emph{%
\ be another non-empty set. We define }$n$\emph{-\textbf{player coopetitive
gain game over the strategy support }}$(E,C)$\emph{\ any pair }
\[
G=(f,>),
\]
\emph{where }$f$\emph{\ is a vector function from the Cartesian product }$%
^{\times }E\times C$\emph{\ (here }$^{\times }E$\emph{\ denotes the classic
strategy-profile space of }$n$\emph{-player normal form games, i.e. the
Cartesian product of the family }$E$\emph{) into the }$n$\emph{-dimensional
Euclidean space }$\Bbb{R}^{n}$\emph{\ and }$>$\emph{\ is the natural
sup-order of this last Euclidean space. The element of the set }$C$\emph{\
will be called \textbf{cooperative strategies of the game}.}

\bigskip

A particular aspect of our coopetitive game model is that any coopetitive game $G$ determines univocally a family of classic normal-form games and vice versa; so that any coopetitive game could be defined as a family of normal-form games. In what follows we precise this very important aspect of the model.

\bigskip

\textbf{Definition (the family of normal-form games associated with a
coopetitive game).}\emph{\ Let }$G=(f,>)$\emph{\ be a coopetitive game over
a strategic support }$(E,C)$\emph{. And let }
\[
g=(g_{z})_{z\in C}
\]
\emph{be the family of classic normal-form games whose member }$g_{z}$\emph{\ is, for any cooperative strategy }$z$\emph{\ in }$C$\emph{, the normal-form game} 
\[
G_{z}:=(f(.,z),>), 
\]
\emph{where the payoff function }$f(.,z)$\emph{\ is the section } 
\[
f(.,z):\;^{\times }E\rightarrow \Bbb{R}^{n} 
\]
\emph{of the function $f$, defined (as usual) by} 
\[
f(.,z)(x)=f(x,z), 
\]
\emph{for every point }$x$\emph{\ in the strategy profile space }$^{\times
}E $\emph{. We call the family }$g$\emph{\ (so defined) \textbf{family of
normal-form games associated with (or determined by) the game} }$G$\emph{ and we call \textbf{normal section} of the game $G$ any member of the family $g$.}

\bigskip

We can prove this (obvious) theorem.

\bigskip

\textbf{Theorem.}\emph{\ The family }$g$\emph{\ of normal-form games
associated with a coopetitive game }$G$\emph{\ uniquely determines the game.
In more rigorous and complete terms, the correspondence }$G\mapsto g$\emph{\
is a bijection of the space of all coopetitive games - over the strategy
support }$(E,C)$\emph{\ - onto the space of all families of normal form
games - over the strategy support }$E$\emph{\ - indexed by the set }$C$\emph{%
.}

\emph{\bigskip }

\emph{Proof.} This depends totally from the fact that we have the following
natural bijection between function spaces: 
\[
\mathcal{F}(^{\times }E\times C,\Bbb{R}^{n})\rightarrow \mathcal{F}(C,%
\mathcal{F}(^{\times }E,\Bbb{R}^{n})):f\mapsto (f(.,z))_{z\in C}, 
\]
which is a classic result of theory of sets. $\blacksquare $

\bigskip

Thus, the exam of a coopetitive game should be equivalent to the exam of a
whole family of normal-form games (in some sense we shall specify).

\bigskip

In this paper we suggest how this latter examination can be conducted and
what are the solutions corresponding to the main concepts of solution which
are known in the literature for the classic normal-form games, in the case of
two-player coopetitive games.

\bigskip

\subsection{\textbf{Two players coopetitive games}}

\bigskip

In this section we specify the definition and related concepts of two-player
coopetitive games; sometimes (for completeness) we shall repeat some
definitions of the preceding section.

\bigskip

\textbf{Definition (of coopetitive game).}\emph{\ Let }$E$\emph{, }$F$\emph{%
\ and }$C$\emph{\ be three nonempty sets. We define \textbf{two player
coopetitive gain game carried by the strategic triple} }$(E,F,C)$\emph{\ any
pair of the form} 
\[
G=(f,>), 
\]
\emph{where }$f$\emph{\ is a function from the Cartesian product }$E\times
F\times C$\emph{\ into the real Euclidean plane }$\Bbb{R}^{2}$\emph{\ and
the binary relation }$>$\emph{\ is the usual sup-order of the Cartesian
plane (defined component-wise, for every couple of points }$p$\emph{\ and }$q
$\emph{, by }$p>q$\emph{\ iff }$p_{i}>q_{i}$\emph{, for each index }$i$\emph{%
).}

\bigskip

\textbf{Remark (coopetitive games and normal form games).} The difference
among a two-player normal-form (gain) game and a two player coopetitive (gain) game
is the fundamental presence of the third strategy Cartesian-factor $C$. The
presence of this third set $C$ determines a total change of perspective with
respect to the usual exam of two-player normal form games, since we now have to consider
a normal form game $G(z)$, for every element $z$ of the set $C$; we have,
then, to study an entire ordered family of normal form games in its own
totality, and we have to define a new manner to study these kind of game
families.

\bigskip

\subsection{\textbf{Terminology and notation}}

\bigskip

\textbf{Definitions.}\emph{\ Let }$G=(f,>)$\emph{\ be a two player
coopetitive gain game carried by the strategic triple }$(E,F,C)$\emph{. We
will use the following terminologies:}

\begin{itemize}
\item  \emph{the function }$f$\emph{\ is called the \textbf{payoff function
of the game} }$G$\emph{;}

\item  \emph{the first component }$f_{1}$\emph{\ of the payoff function }$f$%
\emph{\ is called \textbf{payoff function of the first player} and
analogously the second component }$f_{2}$\emph{\ is called \textbf{payoff
function of the second player};}

\item  \emph{the set }$E$\emph{\ is said \textbf{strategy set of the first
player} and the set }$F$\emph{\ the \textbf{strategy set of the second player%
};}

\item  \emph{the set }$C$\emph{\ is said the \textbf{cooperative (or common)
strategy set of the two players};}

\item  \emph{the Cartesian product }$E\times F\times C$\emph{\ is called the
(\textbf{coopetitive) strategy space of the game} }$G$\emph{.}
\end{itemize}

\bigskip

\textbf{Memento.} The first component $f_{1}$ of the payoff function $f$ of
a coopetitive game $G$ is the function of the strategy space $E\times
F\times C$ of the game $G$ into the real line $\Bbb{R}$ defined by the first
projection 
\[
f_{1}(x,y,z):=\mathrm{pr}_{1}(f(x,y,z)), 
\]
for every strategic triple $(x,y,z)$ in $E\times F\times C$; in a similar fashion we
proceed for the second component $f_{2}$ of the function $f$.

\bigskip

\textbf{Interpretation.} We have:

\begin{itemize}
\item  two players, or better an ordered pair $(1,2)$ of players;

\item  anyone of the two players has a strategy set in which to choose
freely his own strategy;

\item  the two players can/should \emph{cooperatively} choose strategies $z$
in a third common strategy set $C$;

\item  the two players will choose (after the exam of the entire game $G$)
their cooperative strategy $z$ in order to maximize (in some sense we shall
define) the vector gain function $f$.
\end{itemize}

\bigskip

\subsection{\textbf{Normal form games of a coopetitive game}}

\bigskip

Let $G$ be a coopetitive game in the sense of above definitions. For any
cooperative strategy $z$ selected in the cooperative strategy space $C$,
there is a corresponding normal form gain game 
\[
G_{z}=(p(z),>), 
\]
upon the strategy pair $(E,F)$, where the payoff function $p(z)$ is the
section 
\[
f(.,z):E\times F\rightarrow \Bbb{R}^{2}, 
\]
of the payoff function $f$ of the coopetitive game - the section is defined,
as usual, on the competitive strategy space $E\times F$, by 
\[
f(.,z)(x,y)=f(x,y,z), 
\]
for every bi-strategy $(x,y)$ in the bi-strategy space $E\times F$.

\bigskip

Let us formalize the concept of game-family associated with a coopetitive
game.

\bigskip

\textbf{Definition (the family associated with a coopetitive game).}\emph{\
Let }$G=(f,>)$\emph{\ be a two player coopetitive gain game carried by the
strategic triple }$(E,F,C)$\emph{. We naturally can associate with the game }%
$G$\emph{\ a family }$g=(g_{z})_{z\in C}$ \emph{of normal-form games defined
by 
\[
g_{z}:=G_{z}=(f(.,z),>), 
\]
for every }$z$ \emph{in }$C$\emph{, which we shall call \textbf{the family
of normal-form games associated with the coopetitive game} }$G$\emph{.}

\bigskip

\textbf{Remark.} It is clear that with any above family of normal form games 
\[
g=(g_{z})_{z\in C}, 
\]
with $g_{z}=(f(.,z),>)$, we can associate:

\begin{itemize}
\item  a family of payoff spaces 
\[
(\mathrm{im}f(.,z))_{z\in C}, 
\]
with members in the payoff universe $\Bbb{R}^{2}$;

\item  a family of Pareto maximal boundary 
\[
(\partial ^{*}G_{z})_{z\in C}, 
\]
with members contained in the payoff universe $\Bbb{R}^{2}$;

\item  a family of suprema 
\[
(\mathrm{sup}G_{z})_{z\in C}, 
\]
with members belonging to the payoff universe $\Bbb{R}^{2}$;

\item  a family of Nash zones 
\[
(\mathcal{N}(G_{z}))_{z\in C}; 
\]
with members contained in the strategy space $E\times F$;

\item  a family of conservative bi-values 
\[
v^{\#}=(v_{z}^{\#})_{z\in C}; 
\]
in the payoff universe $\Bbb{R}^{2}$.
\end{itemize}

\bigskip

And so on, for every meaningful known feature of a normal form game.

\bigskip

Moreover, we can interpret any of the above families as \emph{set-valued paths} in
the strategy space $E\times F$ or in the payoff universe $\Bbb{R}^{2}$.

\bigskip

It is just the study of these induced families which becomes of great
interest in the examination of a coopetitive game $G$ and which will enable
us to define (or suggest) the various possible solutions of a coopetitive game.

\bigskip

\section{\textbf{Solutions of a coopetitive game}}

\bigskip

\subsection{\textbf{Introduction}}

\bigskip

The two players of a coopetitive game $G$ should choose the cooperative
strategy $z$ in $C$ in order that:

\begin{itemize}
\item  the reasonable Nash equilibria of the game $G_{z}$ are $f$-preferable
than the reasonable Nash equilibria in each other game $G_{z^{\prime }}$;

\item  the supremum of $G_{z}$ is greater (in the sense of the usual order
of the Cartesian plane) than the supremum of any other game $G_{z^{\prime }}$%
;

\item  the Pareto maximal boundary of $G_{z}$ is higher than that of any
other game $G_{z^{\prime }}$;

\item  the Nash bargaining solutions in $G_{z}$ are $f$-preferable than
those in $G_{z^{\prime }}$;

\item  \emph{in general, fixed a common kind of solution for any game }$%
G_{z} $\emph{, say }$S(z)$\emph{\ the set of these kind of solutions for the
game }$G_{z}$\emph{, we can consider the problem to find all the optimal
solutions (in the sense of Pareto) of the set valued path }$S$\emph{,
defined on the cooperative strategy set }$C$\emph{. Then, we should face the
problem of \textbf{selection of reasonable Pareto strategies} in the
set-valued path }$S $\emph{\ via proper selection methods (Nash-bargaining,
Kalai-Smorodinsky and so on).}
\end{itemize}

\bigskip

Moreover, we shall consider the maximal Pareto boundary of the payoff space $%
\mathrm{im}(f)$ as an appropriate zone for the bargaining solutions.

\bigskip

The payoff function of a two person coopetitive game is (as in the case of
normal-form game) a vector valued function with values belonging to the
Cartesian plane $\Bbb{R}^{2}$. We note that in general the above criteria
are multi-criteria and so they will generate multi-criteria optimization
problems.

\bigskip

In this section we shall define rigorously some kind of solution, for two
player coopetitive games, based on a bargaining method, namely a
Kalai-Smorodinsky bargaining type. Hence, first of all, we have to precise
what kind of bargaining method we are going to use.

\bigskip

\subsection{\textbf{Bargaining problems}}

\bigskip

In this paper we shall use the following definition of bargaining problem.

\bigskip

\textbf{Definition (of bargaining problem).}\emph{\ Let }$S$\emph{\ be a
subset of the Cartesian plane }$\Bbb{R}^{2}$\emph{\ and let }$a$\emph{\ and }%
$b$\emph{\ be two points of the plane with the following properties:}

\begin{itemize}
\item  \emph{they belong to the small interval containing }$S$\emph{, if this interval is defined (indeed, it is well defined if and only if $S$ is bounded and it is precisely the interval $[\inf S,\sup S]_{^{\leq }}$);}

\item  \emph{they are such that }$a<b$\emph{;}

\item  \emph{the intersection} 
\[
 [a,b]_{^{\leq }}\cap \partial ^{*}S, 
\]
\emph{among the interval }$[a,b]_{^{\leq }}$\emph{\ with end points }$a$%
\emph{\ and }$b$\emph{\ (it is the set of points greater than }$a$\emph{\
and less than }$b$\emph{, \textbf{it is not} the segment }$[a,b]$\emph{) and
the maximal boundary of }$S$ \emph{is non-empty.}
\end{itemize}

\emph{In this conditions, we call\textbf{\ bargaining problem on} }$S$\emph{%
\ \textbf{corresponding to the pair of extreme points} }$(a,b)$\emph{, the
pair} 
\[
P=(S,(a,b)). 
\]
\emph{Every point in the intersection among the interval }$[a,b]_{^{\leq }}$%
\emph{\ and the Pareto maximal boundary of }$S$\emph{\ is called \textbf{%
possible solution of the problem} }$P$\emph{. Some time the first extreme
point of a bargaining problem is called \textbf{the initial point of the problem}
(or \textbf{disagreement point} or \textbf{threat point}) and the second
extreme point of a bargaining problem is called \textbf{utopia point} of the problem.}

\bigskip

In the above conditions, when $S$ is convex, the problem $P$ is said to be
convex and for this case we can find in the literature many existence
results for solutions of $P$ enjoying prescribed properties
(Kalai-Smorodinsky solutions, Nash bargaining solutions and so on ...).

\bigskip

\textbf{Remark.} Let $S$ be a subset of the Cartesian plane $\Bbb{R}^{2}$
and let $a$ and $b$ two points of the plane belonging to the smallest
interval containing $S$ and such that $a\leq b$. Assume the Pareto maximal
boundary of $S$ be non-empty. If $a$ and $b$ are a lower bound and an upper
bound of the maximal Pareto boundary, respectively, then the intersection 
\[
\lbrack a,b]_{^{\leq }}\cap \partial ^{*}S 
\]
is obviously not empty. In particular, if $a$ and $b$ are the extrema of $S$
(or the extrema of the Pareto boundary $S^{*}=\partial ^{*}S$) we can
consider the following bargaining problem 
\[
P=(S,(a,b)),\;(or\;P=(S^{*},(a,b))) 
\]
and we call this particular problem a \emph{standard bargaining problem on }$%
S$ (or \emph{standard bargaining problem on the Pareto maximal boundary} $%
S^{*}$).

\bigskip

\subsection{\textbf{Kalai solution for bargaining problems}}

\bigskip

Note the following property.

\bigskip

\textbf{Property.}\emph{\ If }$(S,(a,b))$\emph{\ is a bargaining problem
with }$a<b$\emph{, then there is at most one point in the intersection } 
\[
\lbrack a,b]\cap \partial ^{*}S, 
\]
\emph{where }$[a,b]$\emph{\ is the \textbf{segment joining the two points} }$%
a$\emph{\ and }$b$\emph{.}

\emph{\bigskip }

\emph{Proof.} Since if a point $p$ of the segment $[a,b]$ belongs to the
Pareto boundary $\partial ^{*}S$, no other point of the segment itself can
belong to Pareto boundary, since the segment is a totally ordered subset of
the plane (remember that $a<b$). $\blacksquare $

\bigskip

\textbf{Definition (Kalai-Smorodinsky). }\emph{We call \textbf{%
Kalai-Smorodinsky solution} (or \textbf{best compromise solution}) \textbf{%
of the bargaining problem} }$(S,(a,b))$ \emph{the unique point of the
intersection} 
\[
\lbrack a,b]\cap \partial ^{*}S, 
\]
\emph{if this intersection is non empty.}

\bigskip

So, in the above conditions, the Kalai-Smorodinsky solution $k$ (if it
exists) enjoys the following property: there is a real $r$ in $[0,1]$ such
that 
\[
k=a+r(b-a), 
\]
or 
\[
k-a=r(b-a), 
\]
hence 
\[
\frac{k_{2}-a_{2}}{k_{1}-a_{1}}=\frac{b_{2}-a_{2}}{b_{1}-a_{1}}, 
\]
if the above ratios are defined; these last equality is the \emph{%
characteristic property of Kalai-Smorodinsky solutions}.

\bigskip

We end the subsection with the following definition.

\bigskip

\textbf{Definition (of Pareto boundary). }\emph{We call \textbf{Pareto
boundary} every subset }$M$\emph{\ of an ordered space which has only
pairwise incomparable elements.}

\bigskip

\subsection{\textbf{Nash (proper) solution of a coopetitive game}}

\bigskip

Let $N:=\mathcal{N}(G)$ be the union of the Nash-zone family of a
coopetitive game $G$, that is the union of the family $(\mathcal{N}%
(G_{z}))_{z\in C}$ of all Nash-zones of the game family $g=(g_{z})_{z\in C}$
associated to the coopetitive game $G$. We call \emph{Nash path of the game} 
$G$ the multi-valued path 
\[
z\mapsto \mathcal{N}(G_{z}) 
\]
and Nash zone of $G$ the trajectory $N$ of the above multi-path. Let $N^{*}$
be the Pareto maximal boundary of the Nash zone $N$. We can consider the
bargaining problem 
\[
P_{\mathcal{N}}=(N^{*},\inf (N^{*}),\sup (N^{*})). 
\]

\bigskip

\textbf{Definition.}\emph{\ If the above bargaining problem }$P_{\mathcal{N}%
} $\emph{\ has a Kalai-Smorodinsky solution }$k$\emph{, we say that }$k$%
\emph{\ is\textbf{\ the properly coopetitive solution of the coopetitive game%
} }$G$\emph{.}

\bigskip

The term ``properly coopetitive'' is clear:

\begin{itemize}
\item  \emph{this solution }$k$ \emph{is determined by cooperation on the
common strategy set }$C$\emph{\ and to be selfish (competitive in the Nash
sense) on the bi-strategy space }$E\times F$\emph{.}
\end{itemize}

\bigskip

\subsection{\textbf{Bargaining solutions of a coopetitive game}}

\bigskip

It is possible, for coopetitive games, to define other kind of solutions,
which are not properly coopetitive, but realistic and sometime affordable.
These kind of solutions are, we can say, \emph{super-cooperative}.

\bigskip

Let us show some of these kind of solutions.

\bigskip

Consider a coopetitive game $G$ and

\begin{itemize}
\item  its Pareto maximal boundary $M$ and the corresponding pair of extrema 
$(a_{M},b_{M})$;

\item  the Nash zone $\mathcal{N}(G)$ of the game in the payoff space and
its extrema $(a_{N},b_{N})$;

\item  the conservative set-value $G^{\#}$ (the set of all conservative
values of the family $g$ associated with the coopetitive game $G$) and its
extrema $(a^{\#},b^{\#})$.
\end{itemize}

\bigskip

\emph{We call:}

\begin{itemize}
\item  \emph{\textbf{Pareto compromise solution of the game} }$G$\emph{\ the
best compromise solution (K-S solution) of the problem } 
\[
(M,(a_{M},b_{M})), 
\]
\emph{if this solution exists;}

\item  \emph{\textbf{Nash-Pareto compromise solution of the game} }$G$\emph{%
\ the best compromise solution of the problem } 
\[
(M,(b_{N},b_{M})) 
\]
\emph{if this solution exists;}

\item  \emph{\textbf{conservative-Pareto compromise solution of the game} }$%
G $\emph{\ the best compromise of the problem} 
\[
(M,(b^{\#},b_{M})) 
\]
\emph{if this solution exists.}
\end{itemize}

\bigskip

\subsection{\textbf{Transferable utility solutions}}

\bigskip

Other possible compromises we suggest are the following.

\bigskip

Consider \emph{the transferable utility Pareto boundary} $M$ \emph{of the
coopetitive game} $G$, that is the set of all points $p$ in the Euclidean
plane (universe of payoffs), between the extrema of $G$, such that their
sum 
\[
+(p):=p_{1}+p_{2} 
\]
is equal to the maximum value of the addition $+$ of the real line $\Bbb{R}$ over the
payoff space $f(E\times F\times C)$ of the game $G$.

\bigskip

\textbf{Definition (TU Pareto solution).} \emph{We call \textbf{transferable
utility compromise solution of the coopetitive game} }$G$\emph{\ the
solution of any bargaining problem }$(M,(a,b))$\emph{, where}

\begin{itemize}
\item  $a$\emph{\ and }$b$ \emph{are points of the smallest interval
containing the payoff space of }$G$

\item  $b$\emph{\ is a point strongly greater than }$a$\emph{;}

\item  $M$\emph{\ is the transferable utility Pareto boundary of the game }$%
G $\emph{;}

\item  \emph{the points }$a$\emph{\ and }$b$\emph{\ belong to different
half-planes determined by }$M$\emph{.}
\end{itemize}

\bigskip

Note that the above forth axiom is equivalent to require that the segment
joining the points $a$ and $b$ intersect $M$.

\bigskip

\subsection{\textbf{Win-win solutions}}

\bigskip

In the applications, if the game $G$ has a member $G_{0}$ of its family
which can be considered as an ``initial game'' - in the sense that the
pre-coopetitive situation is represented by this normal form game $G_{0}$ -
the aims of our study (following the standard ideas on coopetitive
interactions) are

\begin{itemize}
\item  to ``enlarge the pie'';

\item  to obtain a win-win solution with respect to the initial situation.
\end{itemize}

\bigskip

So that we will choose as a threat point $a$ in TU problem $(M,(a,b))$ the
supremum of the initial game $G_{0}$.

\bigskip

\textbf{Definition (of win-win solution).}\emph{\ Let }$(G,z_{0})$\emph{\ be
a \textbf{coopetitive game with an initial point}, that is a coopetitive
game }$G$\emph{\ with a fixed common strategy }$z_{0}$\emph{\ (of its common
strategy set }$C$\emph{). We call the game }$G_{z_{0}}$\emph{\ as \textbf{%
the initial game of} }$(G,z_{0})$\emph{. We call \textbf{win-win solution of
the game} }$(G,z_{0})$\emph{\ any strategy profile }$s=(x,y,z)$\emph{\ such
that the payoff of }$G$\emph{\ at }$s$\emph{\ is strictly greater than the
supremum }$L$ \emph{of the \textbf{payoff core} of the initial game }$%
G(z_{0})$\emph{.}

\bigskip

\textbf{Remark.} The payoff core of a normal form gain game $G$ is the portion of the Pareto maximal boundary $G^{*}$ of the game which is greater than the conservative bi-value of $G$. 

\bigskip

\textbf{Remark.} From an applicative point of view, the above requirement
(to be strictly greater than $L$) is very strong. More realistically, we can
consider as win-win solutions those strategy profiles which are strictly
greater than any reasonable solution of the initial game $G_{z_{0}}$.

\bigskip

\textbf{Remark.} In particular, observe that, if the collective payoff
function 
\[
^{+}(f)=f_{1}+f_{2} 
\]
has a maximum (on the strategy profile space $S$) strictly greater than the
collective payoff $L_{1}+L_{2}$ at the supremum $L$ of the payoff core of
the game $G_{z_{0}}$, the portion $M(>L)$ of TU Pareto boundary $M$ which is greater than $L$ is
non-void and it is a segment. So that we can choose as a threat point $a$ in
our problem $(M,(a,b))$ the supremum $L$ of the payoff core of the initial
game $G_{0}$ \emph{to obtain some compromise solution}.

\bigskip

\textbf{Standard win-win solution.} A natural choice for the utopia point $b$
is the supremum of the portion $M_{\geq a}$ of the transferable utility
Pareto boundary $M$ which is upon (greater than) this point $a$: 
\[
M_{\geq a}=\{m\in M:m\geq a\}. 
\]

\bigskip

\textbf{Non standard win-win solution.} Another kind of\textbf{\ }solution
can be obtained by choosing $b$ as the supremum of the portion of $M$ that
is bounded between the minimum and maximum value of that player $i$ that
gains more in the coopetitive interaction, in the sense that 
\[
\max (\mathrm{pr}_{i}(\mathrm{im}f))-\max (\mathrm{pr}_{i}(\mathrm{im}%
f_{0}))>\max (\mathrm{pr}_{3-i}(\mathrm{im}f))-\max (\mathrm{pr}_{3-i}(%
\mathrm{im}f_{0})). 
\]

\bigskip

\section{\textbf{Coopetitive games for Greek crisis}}

\bigskip

Our first hypothesis is that Germany must stimulate the domestic demand and
to re-balance its trade surplus in favor of Greece. The second hypothesis is
that Greece, a country with a declining competitiveness of its products and
a small export share, aims at growth by undertaking innovative investments
and by increasing its exports primarily towards Germany and also towards the
other euro countries.

\bigskip

The coopetitive model that we propose hereunder must be interpreted as a
normative model, in the sense that it shows the more appropriate solutions
of a win-win strategy chosen by considering both competitive and cooperative
behaviors.

\bigskip

The strategy spaces of the two models are:

\begin{itemize}
\item  the strategy set of Germany $E$, set of all possible consumptions of
Germany, in our model, given in conventional monetary unit; we shall assume that the strategies of Germany
directly influence only Germany pay-off;

\item  the strategy set of Greece $F$, set of all possible investments of
Greece, in our model, given in conventional monetary unit; we shall assume that the strategies of Greece
directly influence only Greece pay-off;

\item  a shared strategy set $C$, whose elements are determined together by
the two countries, when they choose their own respective strategies $x$ and $y$, Germany and Greece. Every $z$ in $C$ represents an amount - given in conventional monetary unit - of Greek exports imported by Germany.
\end{itemize}

\bigskip

Therefore, in the two models we assume that Germany and Greece define the
set of coopetitive strategies.

\bigskip

\section{\textbf{The mathematical model}}

\bigskip

\textbf{Main Strategic assumptions.} We assume that:

\begin{itemize}
\item  any real number $x$, belonging to the unit interval $U=[0,1]$, can
represent a consumption of Germany (given in conventional monetary unit);

\item  any real number $y$, in the same unit interval $U$, can represent an
investment of Greece (given in conventional monetary unit);

\item  any real number $z$, again in $U$, can be the amount of Greek exports
which is imported by Germany (given in conventional monetary unit).
\end{itemize}

\bigskip

\subsection{Payoff function of Germany}

\bigskip

We assume that the payoff function of Germany $f_{1}$ is its \emph{gross
domestic demand}:

\begin{itemize}
\item  $f_{1}$ is equal to the private consumption function $c_{1}$ plus the gross
investment function $I_{1}$ plus government spending (that we shall assume
equal $0$) plus export function $X_{1}$ minus import function $M_{1}$, that
is 
\[
f_{1}=c_{1}+I_{1}+X_{1}-M_{1}. 
\]
We assume that:

\item  the private consumption function $c_{1}$ is the first projection of
the strategic Cartesian cube $U^{3}$, 
\[
c_{1}(x,y,z)=x, 
\]
since we assume the private consumption of Germany the first strategic
component of strategy profiles in $U^{3}$;

\item  we assume the gross investment function $I_{1}$ constant on the cube $U^{3}$ and by translation we can suppose $I_{1}$ equal zero;

\item  the export function $X_{1}$ is defined by 
\[
X_{1}(x,y,z)=(1+x)^{-1}, 
\]
for every consumption $x$ of Germany; so we assume that the export function $X_{1}$ is a strictly decreasing function with respect to the first argument;

\item  the import function $M_{1}$ is the third projection of the strategic cube, namely
\[
M_{1}(x,y,z)=z, 
\]
since we assume the import function $M$ depending only upon the cooperative
strategy $z$ of the coopetitive game $G$, our third strategic component of
the strategy profiles in $U^{3}$.
\end{itemize}

\bigskip

\textbf{Recap.} We then assume as payoff function of Germany its gross
domestic demand $f_{1}$, which in our model is equal, at every triple $(x,y,z)$ in the strategic cube $U^{3}$, to the sum of the strategies $x$, $-z$ with the export function $X_{1}$, viewed as a reaction function with
respect to the German domestic consumption (so that $f_{1}$ is the difference of
the first and third projection of the Cartesian product $U^{3}$ plus the
function export function $X_{1}$).

\bigskip

Concluding, the payoff function of Germany is the function $f_{1}$ of the
cube $U^{3}$ into the real line $\Bbb{R}$, defined by 
\[
f_{1}(x,y,z)=x+1/(x+1)-z, 
\]
for every triple $(x,y,z)$ in the cube $U^{3}$; where the reaction function $X_{1}$, defined from the unit interval $U$ into the real line $\Bbb{R}$ by 
\[
X_{1}(x)=1/(x+1), 
\]
for every consumption $x$ of Germany in the interval $U$, is the export
function of Germany mapping the level $x$ of consumption into the level $X_{1}(x)$ of German export corresponding to that consumption level $x$.

$\bigskip $

The function $X_{1}$ is a strictly decreasing function, and only this
monotonicity is the relevant property of $X_{1}$ for our coopetitive model.

\bigskip

\subsection{Payoff function of Greece}

\bigskip

We assume that the payoff function of Greece $f_{2}$ is again its gross
domestic demand - private consumption $c_{2}$ plus gross investment $I_{2}$ plus government
spending (assumed to be $0$) plus exports $X_{2}$ minus imports $M_{2}$), 
\[
f_{2}=c_{2}+I_{2}+X_{2}-M_{2}. 
\]
We assume that:

\begin{itemize}
\item  the function $c_{2}$ is irrelevant in our analysis, since we assume
the private consumption independent from the choice of the strategic triple $(x,y,z)$
in the cube $U^{3}$, in other terms we assume the function $c_{2}$ constant on the
cube $U^{3}$ and by translation we can suppose $c_{2}$ itself equal zero;

\item  the function $I_{2}$ is defined by 
\[
I_{2}(x,y,z)=y+nz, 
\]
for every $(x,y,z)$ in $U^{3}$ (see later for the justification);

\item  the export function $X_{2}$ is the linear function defined by 
\[
X_{2}(x,y,z)=z+my, 
\]
for every $(x,y,z)$ in $U^{3}$ (see later for the justification);

\item  the function $M_{2}$ is irrelevant in our analysis, since we assume
the import function independent from the choice of the triple $(x,y,z)$ in $U^{3}$,
in other terms we assume the import function $M_{2}$ constant on the cube $U^{3}$
and by translation we can suppose the import $M_{2}$ equal zero.
\end{itemize}

\bigskip

So the payoff function of Greece is the linear function $f_{2}$ of the cube $U^{3}$ into the real line $\Bbb{R}$, defined by 
\[
f_{2}(x,y,z)=(y+nz)+(z+my)=(1+m)y+(1+n)z, 
\]
for every pair $(x,y,z)$ in the strategic Cartesian cube $U^{3}$.

\bigskip

We note that the function $f_{2}$ does not depend upon the strategies $x$ in 
$U$ chosen by Germany and that $f_{2}$ is a linear function.

\bigskip

The definition of the functions investment $I_{2}$ and export $X_{2}$ must be studied deeply
and carefully, and are fundamental to find the win-win solution.

\begin{itemize}
\item For every investment strategy $y$ in $U$, the term $my$ represents the quantity (monetary) effect of
the Greek investment $y$ on the Greek exports. In fact, the investments,
specially innovative investments, contribute at improving the
competitiveness of Greek goods, favoring the exports.

\item For every cooperative strategy $z$ in $U$, the term $nz$ is the cross-effect of the
cooperative variable $z$ representing the additive level of investment
required to support the production of the production $z$ itself. 

\item We assume the factors $m$ and $n$ strictly positive.

\end{itemize}

\bigskip

\subsection{Payoff function of the game}

\bigskip

We so have build up a coopetitive gain game with payoff function given by 
\begin{eqnarray*}
f(x,y,z) &=&(x+1/(x+1)-z,(1+m)y+z)= \\
&=&(x+1/(x+1),(1+m)y)+z(-1,1+n)
\end{eqnarray*}
for every $x,y,z$ in $[0,1]$.

\bigskip

\subsection{Study of the game $G=(f,>)$}

\bigskip

Note that, fixed a cooperative strategy $z$ in $U$, the section game $G(z)=(p(z),>)$ with payoff function $p(z)$, defined on the square $U\times U$
by 
\[
p(z)(x,y)=f(x,y,z), 
\]
is the translation of the game $G(0)$ by the ``cooperative'' vector 
\[
v(z)=z(-1,1+n), 
\]
so that we can study the initial game $G(0)$ and then we can translate the
various informations of the game $G(0)$ by the vector $v(z)$.

\bigskip

So, let us consider the initial game $G(0)$. The strategy square $S=U^{2}$
of $G(0)$ has vertices $0_{2}$, $e_{1}$, $1_{2}$ and $e_{2}$, where $0_{2}$
is the origin, $e_{1}$ is the first canonical vector $(1,0)$, $1_{2}$ is the
sum of the two canonical vectors $(1,1)$ and $e_{2}$ is the second canonical
vector $(0,1)$.

\bigskip

\subsection{Topological Boundary of the payoff space of $G_{0}$}

\bigskip

In order to determine the Pareto boundary of the payoff space, we shall use
the technics introduced by D. Carf\`{i} in \cite{ca1}. We have 
\[
p_{0}(x,y)=(x+1/(x+1),(1+m)y), 
\]
for every $x,y$ in $[0,1]$. The transformation of the side $[0,e_{1}]$ is
the trace of the (parametric) curve $c:U\rightarrow \Bbb{R}^{2}$ defined by 
\[
c(x)=f(x,0,0)=(x+1/(x+1),0), 
\]
that is the segment 
\[
\lbrack f(0),f(e_{1})]=[(1,0),(3/2,0)]. 
\]
The transformation of the segment $[0,e_{2}]$ is the trace of the curve $%
c:U\rightarrow \Bbb{R}^{2}$ defined by 
\[
c(y)=f(0,y,0)=(1,(1+m)y), 
\]
that is the segment 
\[
\lbrack f(0),f(e_{2})]=[(1,0),(1,1+m)]. 
\]
The transformation of the segment $[e_{1},1_{2}]$ is the trace of the curve $%
c:U\rightarrow \Bbb{R}^{2}$ defined by 
\[
c(y)=f(1,y,0)=(1+1/2,(1+m)y), 
\]
that is the segment 
\[
\lbrack f(e_{1}),f(1_{2})]=[(3/2,0),(3/2,1+m)]. 
\]

\bigskip

\textbf{Critical zone of }$G(0)$\textbf{.} The Critical zone of the game $%
G(0)$ is empty. Indeed the Jacobian matrix is 
\[
J_{f}(x,y)=\left( 
\begin{array}{cc}
1+(1+x)^{-2} & 0 \\ 
0 & 1+m
\end{array}
\right) , 
\]
which is invertible for every $x,y$ in $U$.

\bigskip

\textbf{Payoff space of the game }$G(0)$\textbf{.} So, the payoff space of
the game $G(0)$ is the transformation of the topological boundary of the
strategic square, that is the rectangle with vertices $f(0,0)$, $f(e_{1})$, $%
f(1,1)$ and $f(e_{2})$.

\bigskip

\textbf{Nash equilibria.} The unique Nash equilibrium is the bistrategy $%
(1,1)$. Indeed, 
\[
1+(1+x)^{-2}>0 
\]
so the function $f_{1}$ is increasing with respect to the first argument and
analogously 
\[
1+m>0 
\]
so that the Nash equilibrium is $(1,1)$.

\bigskip

\subsection{\textbf{The payoff space of the coopetitive game} $G$}

\bigskip

The image of the payoff function $f$, is the union of the family of payoff
spaces 
\[
(\mathrm{im}p_{z})_{z\in C}, 
\]
that is the convex envelope of the union of the image $p_{0}(S)$ ($S$ is the square $U\times U$) and of its
translation by the vector $v(1)$, namely the payoff space $p_{1}(S)$: the image of $f$ is an hexagon with vertices  $f(0,0)$, $f(e_{1})$, $f(1,1)$ and their translations by $v(1)$.

\bigskip

\subsection{\textbf{The Pareto maximal boundary of the payoff space of }$G$}

\bigskip

The Pareto sup-boundary of the coopetitive payoffspace $f(S)$ is the segment $[P^{\prime },Q^{\prime }]$, where $P^{\prime }=f(1,1)$ and 
\[
Q^{\prime }=P^{\prime }+v(1). 
\]

\bigskip

\textbf{Possibility of global growth.} It is important to note that the
absolute slope of the Pareto (coopetitive) boundary is $1+n$. Thus the
collective payoff $f_{1}+f_{2}$ of the game is not constant on the Pareto
boundary and, therefore, the game implies the possibility of a global growth.

\bigskip

\textbf{Trivial bargaining solutions.} The Nash bargaining solution on the
segment $[P^{\prime },Q^{\prime }]$ with respect to the infimum of the
Pareto boundary and the Kalai-Smorodinsky bargaining solution on the segment 
$[P^{\prime },Q^{\prime }]$, with respect to the infimum and the supremum of
the Pareto boundary, coincide with the medium point of the segment $%
[P^{\prime },Q^{\prime }]$.
This solution is not acceptable from Germany point of view, it is collectively better than the supremum of $G_{0}$ but it is disadvantageous for Germany (it suffers a loss!): this solution can be thought as a rebalancing solution but it is not realistically implementable.

\bigskip

\subsection{\textbf{Transferable utility solution}}

\bigskip

In this coopetitive context it is more convenient to adopt a transferable
utility solution, indeed:

\begin{itemize}
\item  the point of maximum collective gain on the whole of the coopetitive payoff space is the point 
\[
Q^{\prime }=(1/2,2+m+n).
\]
\end{itemize}

\bigskip

\subsection{\textbf{Rebalancing win-win best compromise solution}}

\bigskip

Thus we propose a rebalancing win-win kind of coopetitive solution, as it follows
(in the case $m=0$):

\bigskip

\begin{itemize}
\item[1.]  we consider the portion $s$ of transferable utility Pareto
boundary 
\[
M:=(0,5/2+n)+\Bbb{R}(1,-1),
\]
obtained by intersecting $M$ itself with the strip determined (spanned by convexifying) by the
straight lines $e_{2}+\Bbb{R}e_{1}$ and 
\[
(2+n)e_{2}+\Bbb{R}e_{1},
\]
\emph{these are the straight lines of maximum gain for Greece in games }$G(0)
$\emph{\ and }$G$\emph{\ respectively}.
\end{itemize}

\bigskip

\begin{itemize}
\item[2.]  we consider the Kalai-Smorodinsky segment $s^{\prime }$ of
vertices $(3/2,1)$ - supremum of the game $G(0)$ - and the supremum of the
segment $s$.
\end{itemize}

\bigskip

\begin{itemize}
\item[3.]  our best payoff coopetitive compromise is the unique point $K$ in
the intersection of segments $s$ and $s^{\prime }$, that is the best compromise
solution of the bargaining problem 
\[
(s,(\sup G_{0},\sup s)).
\]
\end{itemize}

\bigskip

\subsection{\textbf{Win-win solution}}

\bigskip

This best payoff coopetitive compromise $K$ represents a win-win solution
with respect to the initial supremum $(3/2,1)$. So that, as we repeatedly
said, \emph{also Germany can increase its initial profit from coopetition}.

\bigskip

\textbf{Win-win strategy procedure.} The win-win payoff $K$ can be obtained (by chance)
in a properly coopetitive fashion in the following way:

\begin{itemize}
\item  1) the two players agree on the cooperative strategy $1$ of the common set $C$;

\item  2) the two players implement their respective Nash strategies of game 
$G(1)$; the unique Nash equilibrium of $G(1)$ is the bistrategy $(1,1)$;

\item  3) finally, they share the ``social pie'' 
\[
5/2+n=(f_{1}+f_{2})(1,1,1),
\]
in a \textbf{cooperative fashion} (by contract) according to the decomposition $K$.
\end{itemize}

\bigskip

\section{Conclusions}

\bigskip

\begin{itemize}
\item  The model of coopetitive game, provided in the present contribution,
is essentially a \emph{normative model}.

\item  It has showed some feasible solutions in a cooperative perspective to
the Greek crisis.

\item  Our model of coopetition has pointed out the strategies that could
bring to win-win solutions in a cooperative perspective for Greece and Germany.
\end{itemize}

\bigskip

We have found:

\begin{itemize}
\item  a properly coopetitive solution, which is not convenient for Germany,
that is the Kalai-Smorodinsky bargaining solution on the coopetitive Nash
zone, set of all possible Nash equilibria of the coopetitive interaction.

\item  The final remarkable result of the work consists in the
determination of a \emph{win-win solution} by a \emph{new selection method
on the transferable utility Pareto boundary} of the coopetitive game.
\end{itemize}

\bigskip

The solution offered by our coopetitive model:

\begin{itemize}
\item  aims at ``sharing the pie fairly'';

\item  shows a win-win and rebalancing outcome for the two countries, within a coopetitive growth path
represented by a non-constant sum game (for instance on the Pareto boundary of the entire
coopetitive interaction).
\end{itemize}

\bigskip

\begin{itemize}
\item  Our analytical model allow us to find a ``fair'' amount of Greek
exports which Germany must cooperatively import as well as the optimal investments necessary to
improve the Greek economy in this context, thus contributing to growth and to the stability
of both the Greek and Germany economies.
\end{itemize}

\end{document}